\documentclass[12pt]{article}
\usepackage{epsfig}
\usepackage{array}
\usepackage[latin1]{inputenc} 
\usepackage{amsmath}
\usepackage{amsfonts}
\usepackage{amssymb}
\usepackage{enumerate}
\usepackage{caption}
\usepackage{wrapfig}
\usepackage{amsthm}
 \usepackage{url}

\newtheorem*{teorema*}{Teorema}
\newtheorem*{nota*}{Nota}
\usepackage{rotating}

\usepackage{makeidx,multicol}

\def\R{\mathbb R}

\def\lgem{l\raisebox{.5ex}[0cm][0cm]{.}l}

\begin{document}
\markboth{\sl Julià Cufí, Agustí Reventós}{Cauchy-Riemann equations}
\title{A historical review of the Cauchy-Riemann equations and the Cauchy Theorem}
\author{Julià Cufí, Agustí Reventós}
\date{}

    \maketitle
    \tableofcontents

\begin{abstract}
In this article we take a historical tour through the Cauchy-Riemann equations and their relationship with Cauchy's theorem on the independence with respect to the path of the integral of a holomorphic function. Because of its  importance we do a detailed and updated study of the contributions of d'Alembert and Euler to these topics. We also review the Cauchy works about the passage from the real to imaginary by paying attention to some  arguments he uses  that are not clear enough. At the end we comment briefly the evolution of Green's formula and its relation with the above problems. 
\end{abstract}

\section{Introduction}

The aim of this article is to take a historical tour through the Cauchy-Riemann equations and their relationship with Cauchy's theorem on the independence with respect to the path of the integral of a holomorphic function. Remember that these equations establish the relationships
  $$\frac{\partial M}{\partial y}=- \frac{\partial N}{\partial x}, \qquad \frac{\partial N}{\partial y}= \frac{\partial M }{\partial x},$$
for the components $M,N$  of a complex function.

The first place where the above equations appear is in d'Alembert's work of 1752, {\em Essai d'une nouvelle théorie de la résistance des fluides}, \cite{dalembert}. This article is difficult to read and, as far as we know, it has not been thoroughly studied in the historiographical bibliography. Because of its  importance, in Section \ref{1310} we do a detailed and updated study of the part of this article that refers to the Cauchy-Riemann equations.

Later, in 1757, these equations are found in the works of Euler \cite{Euler1}, \cite{Euler2} and \cite{Euler3}, in which he gives his famous equations that regulate the movement of a fluid. A few years later, in a work of 1777,  \cite{Euler}, Euler shows that the Cauchy-Riemann equations are valid for the (complex) derivative of a function.
In fact, it is the first place where these equations are present in the context of functions of a complex variable.

The Cauchy-Riemann equations appear again in 1814 in the work of Cauchy
{\em Mémoire sur les intégrales definis}, \cite{1814}. He proves them for one function that is the complex derivative of another, and he deduces from this  that two integrals are equal even though they have been taken along different paths with the same endpoints, over the boundary of a rectangle, as we explain in Section $2$.

The problem of path independence, one of the most important in the theory of functions of a complex variable, was taken up by Cauchy in 1825 in the article {\em Mémoire sur les intégrales definis prises entre des limits imaginaires}, \cite{ 1825b}. Some of the arguments in this memory are not clear enough, as most historians have acknowledged. However, they can be replaced by alternative and rigorous reasonings as we highlight in Section \ref{0410}.

Currently the easiest way to prove Cauchy's theorem is by using Green's formula. The first time this formula appears stated in the form
\begin{eqnarray*}\label{0707}\int_{\partial R} Pdx+Qdy=\int_{R} (\frac{\partial Q}{\partial x}-\frac{\partial P}{\partial y})dx\,dy\end{eqnarray*}
where $P(x,y), Q(x,y)$ are regular functions of $x,y$ and $R$ is a region of the plane, it is in Cauchy's paper from 1846, \cite {cauchy1846}. Cauchy presents it without proof\footnote{At the beginning of this paper Cauhy says that he will later publish the proof, but he did not, as Katz comments in \cite{katz}.} and uses it to show that if a form $Pdx+Qdy$ is exact then  its integral along a closed path is zero. He goes no further, although the same argument would allow him to prove his theorem on path invariance for the integral of a holomorphic function, as discussed in Section \ref{1110}.

The reason  why Cauchy did not give
a proof  of Green's formula is possibly due to the fact that it  is a plane version of the previous results of Gauss on the volume of a body and of Ostrograsdky on the divergence theorem, which we analyze in Section \ref{1110}.

  The first published proof of Green's formula is due to Riemann in his 1851 dissertation, {\em Grundlagen für eine allgemeine Theorie der
Functionen einer veränderlichen complexen
Grösse}, \cite{riemann}. In this work Riemann also gives the definition of a holomorphic function (he does not call it so)
  and shows that such a function satisfies the Cauchy-Riemann equations. For completeness we review in Section \ref{0407}
 these two Riemann contributions.

\section{Path independence in a complex integral}\label{0410b}
One of the most important results of the theory of functions of complex variable is the fact that the integrals of a holomorphic function along two paths that start at the same point and end at another same point, have the same value. This result is known as Cauchy's theorem.

Gauss can be considered the first to discuss this subject since in a letter to Bessel on December 18, 1811 he affirms, speaking about integrals between complex limits, that ``the integral $\int \varphi(x) \,dx$ has only one value, even if taken on different paths''. He says that it is a beautiful theorem of easy proof, and announces  that he will give this proof on another occasion, but he did not (\cite{gauss} p.157 and \cite{Kline}, p.632).
Specifically he says:
\begin{quote}
Now what should one think of $\int \varphi(x)\,dx$ for $x=a+bi$? Obviously, if we want to begin  from clear concepts, we must assume that $x$ passes  through infinitely small increments
(each of the form $\alpha+ \beta i$)
from the value for which the integral is $0$, to $x=a+bi$, and then sum all the  $\varphi(x)\,dx$. In this way the meaning is completely established. But the passage can occur in infinitely many ways: just as one can think of the entire domain of all real magnitudes as  an infinite straight line, so one can make the entire domain of all magnitudes, real and imaginary, meaningful as  an infinite plane, wherein each point  determined by abscissa $= a$ and ordinate  $=b$, represents the quantity
$a+bi$ as it were. The continuous  passage from one value of $x$ to another $a+bi$ accordingly occurs along a curve  and is consequently  possible in infinitely many ways.

I now assert that the integral $\int \varphi(x)\,dx$ always maintains a single value after two different passages, if  $\varphi(x)$ nowhere $=\infty$  within the region enclosed between the curves  representing the two passages.

This is a very beautiful  theorem
   for  which I  will give a not difficult proof at a suitable opportunity.\footnote{This English translation  of a portion of Gauss's letter is due to Ruch, \cite{Ruch}.}
   \end{quote}

As we said in the Introduction, Cauchy deals with the problem of path independence in two different stages: first in a paper from the year 1814 in which he proves the Cauchy-Riemann equations and the path independence for the particular case of the bounadary of a rectangle, and later in a second article from 1824 where he considers  the general case.

  The first of these articles, entitled {\em Mémoire sur les intégrales definis}, \cite{1814}, was presented to the Académie in 1814 but was not published until 1827, and in it Cauchy aims to pass from real integrals to imaginary integrals. Specifically on page 330 he says: I hope to establish the passage from the real to the imaginary through direct and rigorous analysis.

  It is in this article, in a footnote on page 340\footnote{Although Cauchy read his memoir in 1814, this footnote was added to the version presented in 1825 in ``Mémoires'', published in 1827 (see Smithies \cite{Smithies}, p.60).}, where Cauchy gives a formula showing that  two integrals are equal,  even though they have been made along different paths with the same ends, specifically following the sides of a rectangle. This result appears again in {\em Résumé des le\c cons sur le calcul infinitésimal}, Lesson 34, p.204, in {\cite{resume} and in \cite{1822b} p.291.
The formula he gives is
\begin{eqnarray}\label{2403}\int_{x_{0}}^{X}\bigg(f(x+Yi)-f(x+y_{0}i)\bigg)\,dx=\int_{y_{0}}^{Y}(f(X+yi)-f(x_{0}+yi)) \,i\,dy,\end{eqnarray}
which actually means
\begin{eqnarray*}\label{0310}\int_{ADC}f(z)dz=\int_{ABC}f(z)dz\end{eqnarray*}
where $A,B,C,D$ are the vertices of a rectangle traveled in a direct sense.

The proof of equality \eqref{2403} is based on the Cauchy-Riemann equations that Cauchy himself has shown, for a function that is the complex derivative of another,
   in the first part of \cite{1814}, p.336,\label{0910} in the section entitled {\em Des équations qui autorisent le passage du réel a l'imaginaire}. He considers a function $f(y)$ and assumes that the variable $y$ is in turn a function of two new variables $x,z$ and he defines the function
$$F(y)=\int f(y)dy$$ so that $F'(y)=f(y)$. Thus,
\begin{eqnarray*}
\frac{\partial F}{\partial x}&=&f(y)\frac{\partial y}{\partial x}\\\\
\frac{\partial F}{\partial z}&=&f(y)\frac{\partial y}{\partial z}.
\end{eqnarray*}
Hence 
$$\frac{\partial }{\partial z}(f(y)\frac{\partial y}{\partial x})=\frac{\partial }{\partial x}(f(y)\frac {\partial y}{\partial z}).$$

In the particular case where $y(x,z)=x+iz$ it results
$$\frac{\partial }{\partial z}(f(y))=\frac{\partial }{\partial x}(f(y)i)$$
and putting  $f(y)=Q(x,z)+i R(x,z)$
one has
$$\frac{\partial }{\partial z}(Q+iR)=\frac{\partial }{\partial x}(iQ-R)$$
and therefore
\begin{eqnarray*}
\frac{\partial Q}{\partial x}&=&\frac{\partial R}{\partial z}\\\\
\frac{\partial Q}{\partial z}&=&-\frac{\partial R}{\partial x},
\end{eqnarray*}
which are the Cauchy-Riemann equations for the components of $f$.

To prove  \eqref{2403}  Cauchy, in  page 338 of \cite{1814}, uses the  equality
$$\int_{y_{0}}^{Y}\int_{x_{0}}^{X}h(x,y)dx\,dy=\int_{x_{0}}^{X}\int_{y_{0}}^{Y}h(x,y)dy\,dx,$$
first applied to a function $h$ such that
$$h=\frac{\partial S}{\partial x}=\frac{\partial T}{\partial y}, $$
and then to a function $h$ such that
$$h=\frac{\partial S}{\partial y}=-\frac{\partial T}{\partial x}.$$
In the first case
one has 
\begin{eqnarray}\label{0206}\int_{y_{0}}^{Y}(S(X,y)-S(x_{0},y))dy=\int_{x_{0}} ^{X}(T(x,Y)-T(x,y_{0}))dx\end{eqnarray}
  and in the second one 
\begin{eqnarray}\label{0206b}\int_{x_{0}}^{X}(S(x,Y)-S(x,y_{0}))dx=-\int_{y_{0} }^{Y}(T(X,y)-T(x_{0}, y))dy.\end{eqnarray}
  Multiplying \eqref{0206} by $i$ and adding \eqref{0206b} one gets 
  \eqref{2403} for $f=S+iT$.
 
  In other words, the equality \eqref{2403} is true for all functions of the type $f=S+iT$ as long as  $S$ and $T$ satisfy
  the equalities $\frac{\partial S}{\partial x}=\frac{\partial T}{\partial y},\quad \frac{\partial S}{\partial y}=-\frac{\partial T}{\partial x}$, which are the Cauchy-Riemann equations, that are satisfied when $f$ is the complex derivative of another function.\label{0407b}

\medskip
As we discussed before, Cauchy comes again to 
the problem on the independence of the path in the integral of a complex function in \cite{1825b} and \cite{1825a}, which is the same paper, presented at the Académie in 1824, and published in two different journals and in different years. In Section \ref{0410} we will discuss more broadly Cauchy's  contribution to this topic.

  Cauchy's work on the integral of a function over the boundary of a rectangle, which we have just commented, is also discussed in Kline, \cite{Kline}, p.636, and in Smithies \cite{Smithies} , p.60.

\section{The Cauchy-Riemann equations}\label{1310}

The first place where the equations known today as the Cauchy-Riemann equations appear is in d'Alembert's work from 1752, {\em Essai d'une nouvelle théorie de la résistance des fluides}, \cite{dalembert}. The functions for which d'Alembert establishes these equations are the horizontal and vertical components of the force acting on a fluid. He also shows that there is a relationship with the functions of a complex variable.

Later on, Euler also found them in fluid theory in  a paper from 1757,\cite{Euler2} and  again in an article on complex integration from 1777, \cite{Euler} (not published until 1797).

In the following we will comment on the contributions of d'Alembert and Euler on the aspects of fluid theory related to  Cauchy-Riemann equations. However we will begin making reference to a  previous work by Clairaut \cite{Clairaut1740} closely related to the subject discussed by these two authors.

In 1740 Clairaut in the article \cite{Clairaut1740} proved that if a differential form $Pdx+Qdy$ is exact, that is, if there exists a function $f$ such that $df=Pdx+Qdy$, then it must be $P_{y}=Q_{x}$. We note that this statement is equivalent to the equality $f_{xy}=f_{yx}$, later known as the Schwarz identity.

Clairaut's reasoning is as follows. Since $P(x,y)=f_{x}(x,y)$ and $Q(x,y)=f_{y}(x,y)$, integrating it results
$$\int P(x,y)dx=f(x,y)+Y(y),\qquad \int Q(x,y)dy=f(x,y)+X(x)$$
for certain functions $X(x), Y(y).$ Therefore \begin{eqnarray*}\int P(x,y)dx-Y(y)=\int Q(x,y)dy-X(x)\end{eqnarray*} and deriving this equality with respect to $y$ one has
$$\int P_{y}(x,y)dx-Y'(y)=Q(x,y)$$ and now deriving with respect to $x$,
$$P_{y}(x,y)=Q_{x}(x,y),$$
as asserted.

He also said, wrongly, that the condition $P_{y}=Q_{x}$ was sufficient so that the form $Pdx+
Qdy$ was exact, but this is not correct as it was pointed out by  d'Alembert in 1768 (\cite{dalembert1768}, p.14), taking

$$P(x,y)=\frac{y}{x^{2}+y^{2}},\quad Q(x,y)=\frac{-x}{x^{2}+ y^{2}}.$$

Later, Clairaut in 1743  (\cite{Clairaut1743}, p.38), studying the figure of the Earth, affirms  that
in order that a fluid would  be in equilibrium 
it is necessary the form $Pdx+Qdy$, where   $(P,Q)$ are the components  of the force acting on the fluid,
be exact. In particular $Q_{x}=P_{y}$ which is the first Cauchy-Riemann equation for the function $Q+iP$.

\subsection{The work of d'Alembert on fluid theory}\label{1206}
In section \S 19 of the work {\em Essai d'une nouvelle théorie de la résistance des fluides}, \cite{dalembert}, d'Alembert meets again  Clairaut's equality  for the horizontal and vertical components of the 
forces acting on the particles of a fluid in equilibrium.


Later, in this same article, d'Alembert considers the case of a fluid that encounters an obstacle of revolution whose axis has the direction of the fluid.
Then he sees that the components of the fluid velocity, $(q(x,z), p(x,z))$ satisfy the equations
\begin{eqnarray}\label{1910}
\frac{\partial p}{\partial x}=\frac{\partial q}{\partial z},\qquad \frac{\partial p}{\partial z}=-\frac{\partial q}{ \partial x}-\frac{p}{z}
\end{eqnarray}
which would be the Cauchy-Riemann equations for the function $p+iq$ if it was  not for the term $-p/z$ appearing  in the second equation.

Now one coud obtain a description of the movement of the fluid by integrating these equations. Since this is complicated, d'Alembert in 
section \S 58 considers the simpler problem of characterizing pairs of functions $M(x,z), N(x,z)$ such that the differentials $M\,dx+N\,dz$ and $N\,dx-M\,dz$
are exact, that is, that there exist\label{2010} functions $p,q$ with $dp=M\,dx+N\,dz$ and $dq=N\,dx-M\,dz$. These equalities
mean that the functions $p,q$ satisfy \eqref{1910} without the term $-p/z$.

The physical conditions that lead to obtain equations like those of \eqref{1910} but without the term $-p/z$ are given
when the obstacle encountered by the fluid is a cylinder whose axis is perpendicular to the direction of the fluid. Although a cylinder is a surface of revolution, this case is not included in the previous situation since now the axis of the cylinder does not have the direction of the fluid. D'Alembert considers this case 
  in section \S 73 of \cite{dalembert}, where he says that the calculations he had done for surfaces of revolution with the axis following the direction of the fluid can be redone to obtain the desired equations, but he does not.\label{1207b} We we will do this from Euler's equations on page \pageref{3110}.


Below we will present d'Alembert's arguments trying to make them as understandable as possible. For this we will update the notation and some aspects of d'Alembert's reasoning, especially the use of derivatives and the mean value theorem.

\subsection*{Components of the forces acting on a fluid}
Let us consider Figure 5 of d'Alembert, where the points $M,N,O,Q$
represent four fluid particles infinitely close to each other. The point $A$
is an arbitrary point that can be inside or outside the fluid. The lines $AP,MO,NQ$ are parallel and the angle $\angle APM$ is a right angle. Denote $AP=x$, $PM=y$.

\begin{center}
\includegraphics[width=.4\textwidth]{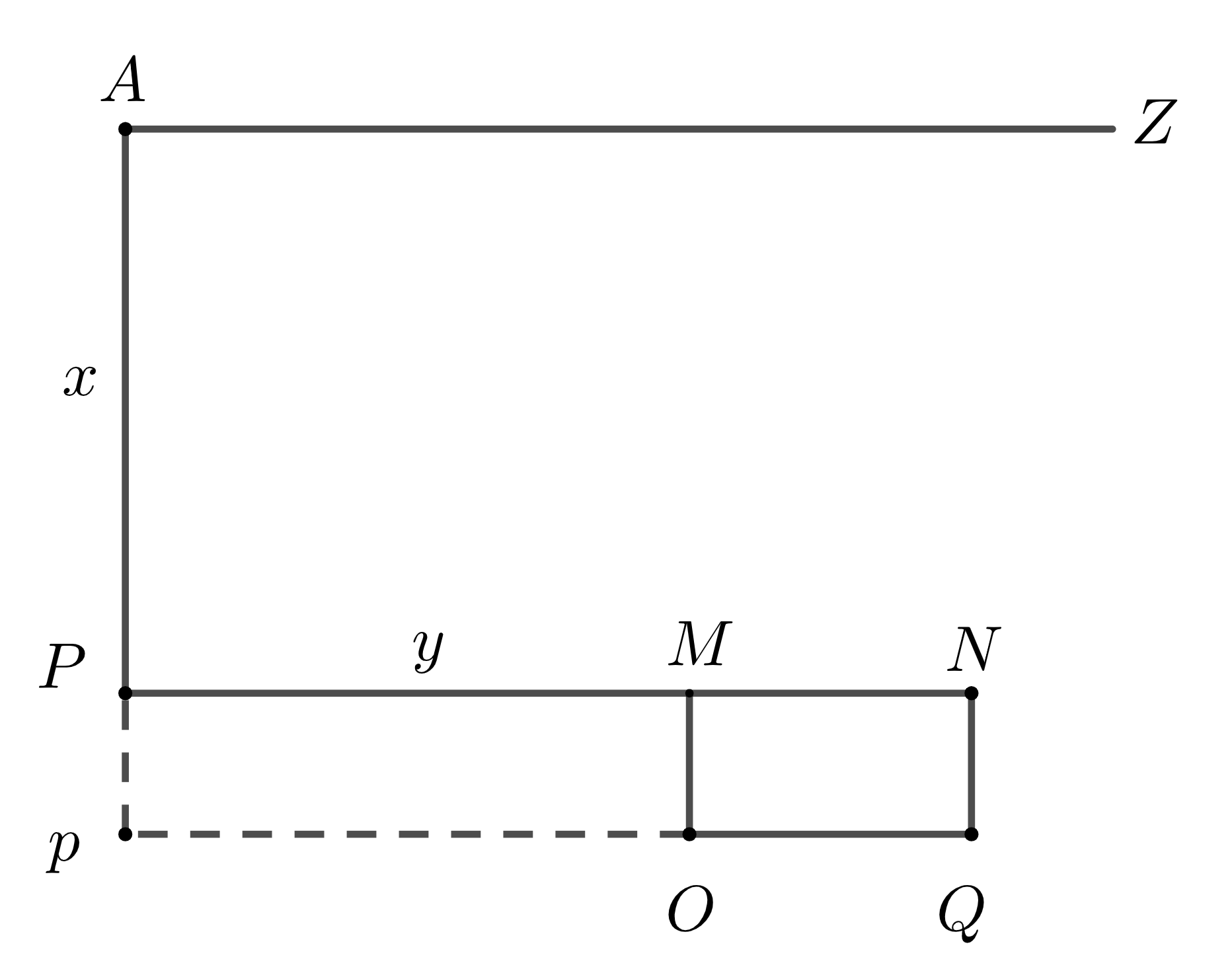}
\centerline{\em Figure 5 of d'Alembert}
\end{center}
Let us break down the forces acting on each of these four points into their horizontal ($PM$ direction) and vertical ($MO$ direction) components.

Denote by $R$ and $Q$ the vertical and horizontal components of the force on $M$, and set $MO=\alpha$, $MN={\cal C}$.

With the current notation, not used by d'Alembert, we can assume that at the point $A$ we have an orthonormal reference $(e_{1},e_{2})$, with $e_{1}$ in the direction of the vector $\overrightarrow{AP}$ and $e_{2}$ in the direction of the vector $\overrightarrow{PM}$. In this way we can write the force vector $\vec{F}$ at the point $M$ like $\vec{F}=Re_{1}+Qe_{2}$, as well as $\overrightarrow{MO}=\alpha e_{1},
\overrightarrow{MN}={\cal C}e_{2}.
$
 
If we think of $\vec{F}$ as a function of the coordinates $x,y$ we can write 
\begin{eqnarray*} \vec{F}(x,y)&=&R(x,y)e_{1}+Q(x,y)e_{2}\quad \mbox{(force at M)}\\
 \vec{F}(x,y+{\cal C})&=&R(x,y+{\cal C})e_{1}+Q( x,y+{\cal C})e_{2}\quad \mbox{(force at  $N$)}\\
\vec{F}(x+\alpha,y)&=&R(x+\alpha,y)e_{1}+Q(x+\alpha,y) e_{2}\quad \mbox{(force at $O$)}\\
\vec{F}(x+\alpha,y+{\cal C})&=&R(x+\alpha,y+{\cal C})e_{1 }+Q(x+\alpha,y+{\cal C})e_{2} \quad \mbox{(force at $Q$)}.
\end{eqnarray*}
By the mean value theorem we have the following expressions for the forces acting on the channels $[MN]$ and $[MO]$.
 
\begin{eqnarray*}
\,[MN]:\quad R(x,y+{\cal C})-R(x,y)&=&{\cal C}\frac{\partial R}{\partial y}(x, y+\xi)\qquad 0\leq \xi\leq {\cal C}\\
\,[MO]:\quad Q(x+\alpha,y)-Q(x,y)&=&\alpha\frac{\partial Q}{\partial x}(x+\eta,y),\qquad 0\leq \eta\leq \alpha,
\end{eqnarray*}
equalities that can be read, as d'Alembert does, saying:

\medskip
$[MN]$: The force at the point $N$ following $NQ$ (the one denoted by $R(x,y+{\cal C})$) will be $ R+{\cal C}\frac{\partial R}{\partial y}$. We have written  $R=R(x,y)$ so as not to overload the notation.

\medskip
$[MO]$: The force at the point $O$ following $OQ$ (the one denoted by $Q(x+\alpha,y)$)
will be $Q+\alpha\frac{\partial Q}{\partial x}.$

\medskip
If the fluid is not homogeneous, these forces must be multiplied by the density $\delta$ which is a function of the point. However, to simplify the calculations, we will assume from now on that $\delta$ is constant, since it is in this situation that the first of the Cauchy-Riemann equations for the functions $R$ and $Q$ appears .

\medskip

By the fundamental principle of fluids in equilibrium, the total sum of forces on the channel $MNOQ$ must vanish. The inner particles cancel their forces with each other so that what must happen is that the sum of forces at the boundary is zero. Equivalently, as d'Alembert says, the strength of the columns $MN$ and $NQ$, following $MN$ and $NQ$ respectively, must be equal to the strength of the columns $MO$ and $ OQ$, following $MO$ and $OQ$ respectively.

The force of the column $MN$ following $MN$ (that is, the second component of this force) is equal to
$$\int_{0}^{{\cal C}}Q(x,y+t)\,dt={\cal C}Q(x,y+\xi).$$

The force of the column $NQ$ following $NQ$ (that is, the first component of the force) is equal to

$$\int_{0}^{\alpha}\bigg(R(x+t,y)+{\cal C}\frac{\partial R}{\partial y}(x+t, y+\xi( t))\bigg)dt=\alpha \bigg (R(x+\eta,y)+{\cal C}\frac{\partial R}{\partial y}(x+\eta,y+\xi(\eta ))\big).$$

Similarly the force of $OQ$ following $OQ$ is equal to
$$\int_{0}^{{\cal C}}\bigg(Q(x,y+t)+\alpha\frac{\partial Q}{\partial x}(x+\eta(t),y +t)\bigg)dt={\cal C}\bigg(Q(x,y+\xi)+\alpha\frac{\partial Q}{\partial x}(x+\eta(\xi),y+\ xi)\big).$$

Finally, the force  of the column $MO$ following $MO$ is equal to $$\alpha R(x+\eta,y).$$

We note that $y\leq y+\xi, y+\bar{\xi}\leq y+{\cal C}$, $x\leq x+\eta, x+\bar{\eta\leq x+\alpha}$.

\medskip

Equalizing the sum of the forces of the columns
$MN$ and $NQ$ with the sum of the forces of the columns $MO$ and $OQ$, we have

\begin{eqnarray*}&&\alpha R(x+\eta,y)+{\cal C}Q(x,y+\xi)+\alpha {\cal C} \frac{\partial Q}{\partial x }(x+\eta(\xi),y+\xi))\\&=&{\cal C}Q(x,y+\bar{\xi})+\alpha R(x+\bar{\eta}, y)+\alpha {\cal C}\frac{\partial R}{\partial y}(x+\bar{\eta},y+\bar{\xi}(\bar{\eta}))).\end{eqnarray*}

Making ${\cal C}$ tend to zero it follows that  $R(x+\eta,y)=R(x+\bar{\eta},y)$ and therefore

\begin{eqnarray*}&&Q(x,y+\xi)+\alpha \frac{\partial Q}{\partial x}(x+\eta(\xi),y+\xi))\\&=&Q(x ,y+\bar{\xi})+\alpha \frac{\partial R}{\partial y}(x+\bar{\eta},y+\bar{\xi}(\bar{\eta}))) .\end{eqnarray*}

Making $\alpha$ tend to zero it follows that  $Q(x,y+\xi)=Q(x,y+\bar{\xi})$ and therefore

\begin{eqnarray*} \frac{\partial Q}{\partial x}(x+\eta(\xi),y+\xi))= \frac{\partial R}{\partial y}(x+\bar{ \eta},y+\bar{\xi}(\bar{\eta}))).\end{eqnarray*}

And finally, if $\alpha $ and ${\cal C}$ tend to zero,  we get

$$\frac{\partial Q}{\partial x}(x,y)=\frac{\partial R}{\partial y}(x,y),$$ which is the first Cauchy-Riemann equation for  the function $Q+iR$.

\subsection*{Pressure, Velocity and Acceleration}
In this section we present the formula for the pressure of a fluid as a function of its speed and an expression for the acceleration that we will need later.

In paragraph \S 27 of the quoted article \cite{dalembert}, d'Alembert considers a  section of a uniform cylindrical channel that widens from the points $A, B$ to the points $P,M$ between  two  walls $AP$ and $BM$ (Figure $10$). He assumes that the fluid is uniform and that in the 
cylindrical part, that is above the section $AB$, the speed is constant and equal to $v_{A}$. In this situation d'Alembert calculates the
velocity, $v_{P}$, and the pressure at any point $P$ in the channel.

\begin{center}
\includegraphics[width=.3\textwidth]{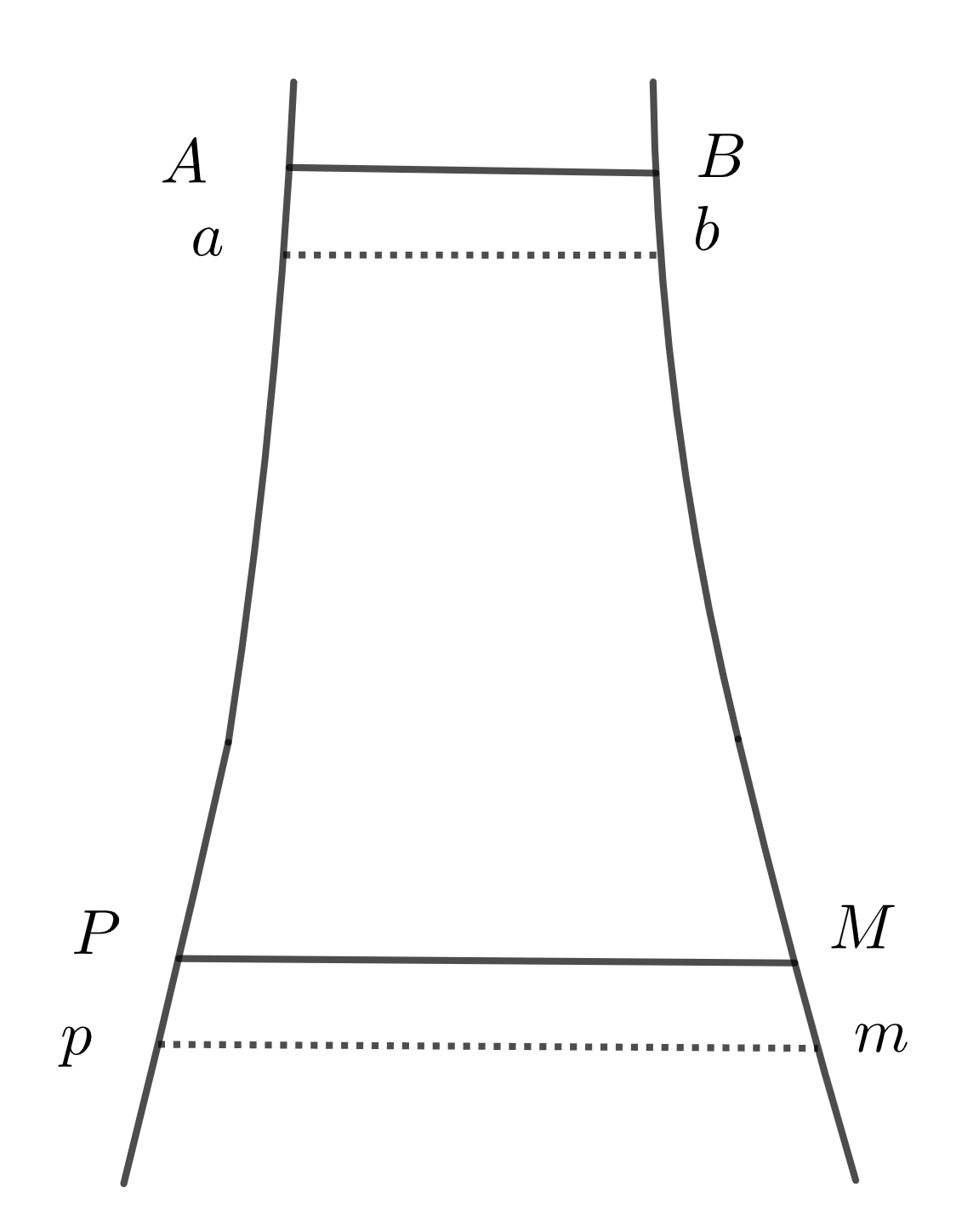}
\centerline{\em Figure 10 of d'Alembert.}
\end{center}

The fluid section $AB$ becomes the 
section $ab$ after an infinitesimal time increment $\Delta t$, and similarly the section $PM$ becomes 
 the section $pm$  in the same time increment.

Equalizing the areas of the infinitesimal rectangles $ABab$ and $PMpm$ one obtains 
$$AB\cdot v_{A}\Delta t=MP\cdot v_{P}\Delta t$$
and therefore,
$$v_{P}=\frac{AB}{PM}v_{A}.$$

To study the pressure in $P$ he  first  observes that the acceleration is negative since when the section of the tube increases the speed decreases.

If $s=s(t)$ measures the length
traveled  by a particle on the channel $AP$
in time $t$, from the point $A$,
  the infinitesimal acceleration at the point $s(t)$ is $-dv/dt$ with $v=ds/dt$,
and the infinitesimal pressure
is $-dv/dt\cdot ds$.

Therefore the pressure at a point $P$ of coordinate $s$ is 
\begin{eqnarray}\label{1611}\int_{0}^{s}-\frac{dv}{dt}\, ds\end{eqnarray}
and hence it is given by
\begin{eqnarray}\label{2905}\int_{0}^{s}-\frac{dv}{dt}\, ds=-\int_{0}^{t}\frac{dv}{dt} v\,dt=\frac{v_{A}^{2}-v(s)^{2}}{2}.\end{eqnarray}

\medskip
In paragraphs \S 43 and \S 44 of the same paper, d'Alembert calculates the horizontal and vertical components of the force acting on a particle $N$ of the fluid. Let us remember that d'Alembert deals with stationary fluids, those in which the velocity depends only on the point and not on the time, that is $v(x,z)=a(q(x,z), p(x ,z))$ for a certain constant $a$. If $(x(t), z(t))$ is a parametrization of the curve described by the particle then  its velocity satisfies
$$(x'(t),z'(t))=a(q(x(t),z(t)), p(x(t),z(t))).$$

Consequently, if we put $dq=Adx+Bdz$, $dp=A'dx+B'dz$ the force that acts on the fluid, which we equate to acceleration, is
\begin{eqnarray}\label{2605b}dv&=&a(dq,dp)=a(Adx+Bdz, A'dx+B'dz)\nonumber\\&=&a^{2}(Aq+Bp, A 'q+B'p).\end{eqnarray}

\subsection*{Fluids encountering an obstacle of revolution}

Next, in paragraph \S 45,  d'Alembert assumes that the fluid encounters an obstacle formed by a body of revolution obtained by rotating the section $ANm$ around the $x$  axis,  that has  the direction of the fluid (Figure 17).
One can think of this section as the graph of a function $z=z(x)$.

\begin{center}
\includegraphics[width=.5\textwidth]{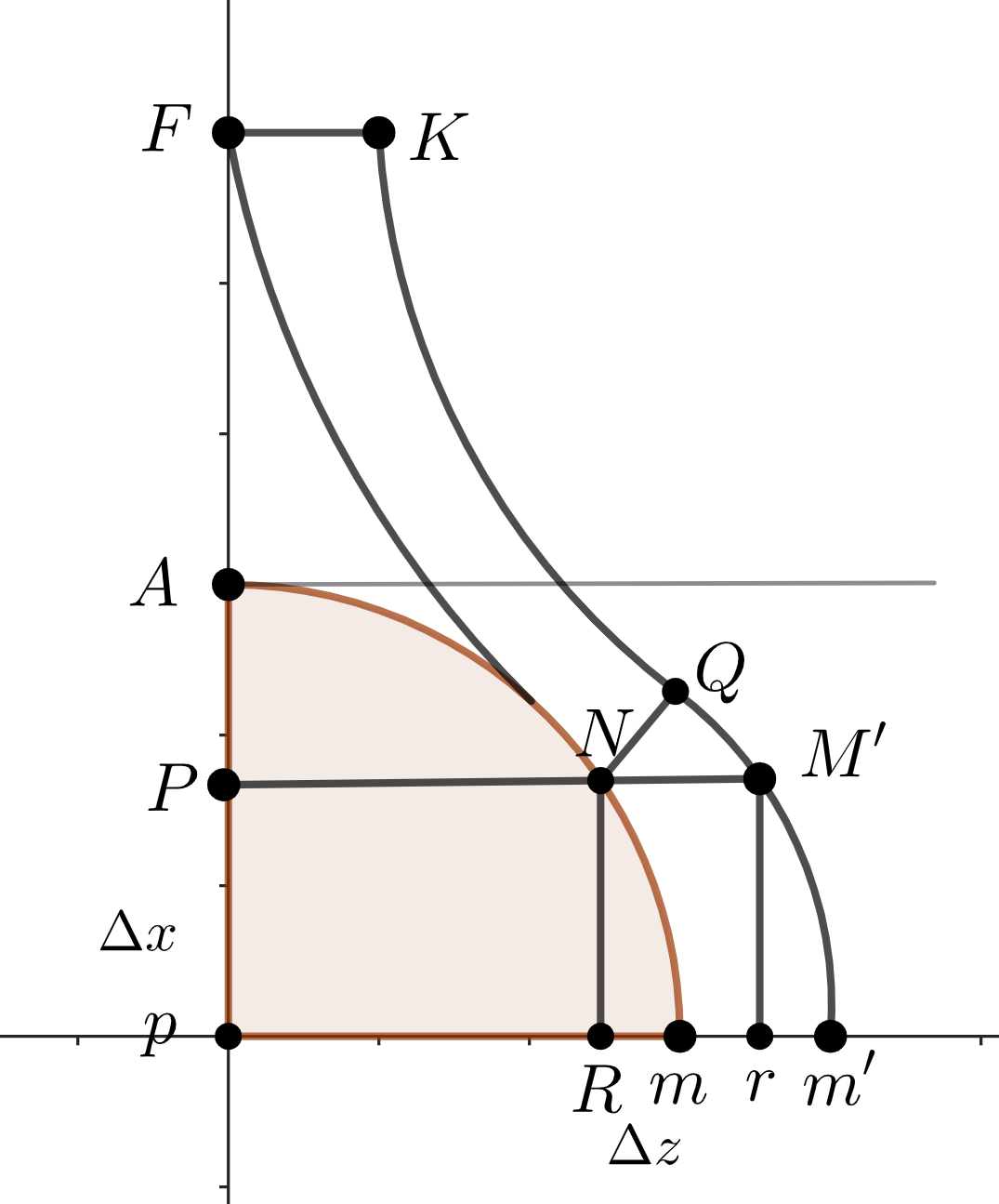}
\centerline{\em Figure 17 of d'Alembert.}
\end{center}

Then he studies  the components of the velocity $v(x,z)=a(q(x,z),p(x,z))$ of the trajectory $FNm$ tangent to the obstacle, and proves, with the notation introduced just before  formula \eqref{2605b}, that $$A'=B, \quad B'=-A-\frac{p}{z}.$$ Since  $A=\frac{\partial q}{\partial x}$, $B=\frac{\partial q}{\partial z}$, $A'=\frac{\partial p}{\partial x}$ and $B'=\frac{\partial p}{\partial z}$, the previous equalities become
\begin{eqnarray}\label{1206b}
\frac{\partial p}{\partial x}&=&\frac{\partial q}{\partial z},\nonumber\\ \frac{\partial p}{\partial z}&=&-\frac {\partial q}{\partial x}-\frac{p}{z}.
\end{eqnarray}

\medskip
To prove the equalities $B'=-A-\frac{p}{z}$ and $A'=B$ d'Alembert uses the principle stating that the volume of fluid that crosses the disc of radius $\rho=FK$, center $F$ and perpendicular to $AP$, in the  time $\Delta t$
is equal to the volume of the fluid that crosses the truncated cone generated by rotating $NQ$ around $AP$ in the same time $\Delta t$.

The first of these two volumes is the volume of a cylinder of radius $\rho$ and height $a\Delta t$, since we are assuming that the velocity far  from the obstacle, has only vertical component of modulus $a$. Therefore, this volume is equal to $\pi \rho^{2}a\Delta t$.

Regarding the second volume, we note that 
it is equal to the volume of a straight cylinder with basis the circular crown of inner radius $z(x)=PN$, width $\delta(x)=NM'$, and height $aq(x,z(x))\Delta t$, since $aq(x,z(x))$ is the vertical component of the velocity at the point $N=(x,z(x))$.
Hence this volume  is $(\pi \delta^{2}+2\pi \delta z)aq\Delta t$.

%
%
%
%

Therefore, the equality of volumes gives $$\pi \rho^{2} a \Delta t=(\pi \delta^{2}+2\pi \delta z)aq\Delta t$$ and consequently
\begin{eqnarray}\label{1610}\delta(x)&=&-z(x)+\sqrt{z(x)^{2}+\frac{\rho^{2}}{q(x ,z(x))}}.\end{eqnarray}

Using Taylor's development at the second order one gets 
\begin{eqnarray}\label{1411}\delta(x)=\frac{\rho^{2}}{2qz}+o(\rho^{2}), \quad \rho\rightarrow 0.\end{eqnarray}

Now d 'Alembert, through
Figure 17, finds the following relation between $\delta'(x)$ and $z'(x)$: 

\begin{eqnarray*}
\delta'(x)=\lim_{\Delta x\to 0}\frac{\delta(x+\Delta x)-\delta(x)}{\Delta x}&=&\lim_{\Delta x\to 0}\frac{mm '-Rr}{\Delta x}\\
&=&\lim_{\Delta x\to 0}\frac{rm'-Rm}{\Delta x}=\lim_{\Delta x\to 0}\frac{rm'-\Delta z}{\Delta x}\\&=&
\lim_{\Delta x\to 0}\frac{rm'}{\Delta x}-z'(x).
\end{eqnarray*}

But since the velocity vector at point $M'$ is tangent to the path $KM'm'$ we have 
$$\lim_{\Delta x\to 0}\frac{rm'}{\Delta x}=\frac{p(x,z(x)+\delta(x))}{q(x,z( x)+\delta(x))}=\frac{p(x,z(x))+\delta(x)p_{z}(x,z(x)+\xi)}{q(x, z(x))+\delta(x)q_{z}(x,z(x)+\eta)}=z'(x)+\delta'(x),$$ with $0<\xi, \eta< \delta(x)$. Since $z'(x)=p(x,z(x))/q(x,z(x))$  we can write, omitting the variables,  the previous equality as

$$\frac{p+\delta p_{z}}{q+\delta q_{z}}=\frac{p}{q}+\delta'$$
and so
\begin{eqnarray}\label{2605}\delta p_{z}=q\delta'+\delta\delta'q_{z}+\frac{\delta p q_{z}}{q}.\end{eqnarray}

%
%

We want to differentiate\footnote{D'Alembert does not use derivatives but arguments with infinitesimals  of the first and second order.}  equality \eqref{2605} twice with respect to $\rho$, at $\rho=0$. 
First of all remark that from \eqref{1411} it comes
\begin{eqnarray*}
\frac{d\delta(x)}{d\rho}=\frac{\rho}{qz}+o(\rho), \qquad 
\frac{d^{2}\delta(x)}{\delta\rho^{2}}_{|\rho=0}=\frac{1}{qz}.
\end{eqnarray*}

We note also that $p_{z}$ and $q_{z}$ depend on $\rho$ through $\xi$ and $\eta$. However, their derivatives with respect to $\rho$ do not play any role since they always appear multiplied by $\delta$ or $\delta'$, which tend to zero with $\rho.$

The second derivative of the term on the left-hand side of equation \eqref{2605} is
\begin{eqnarray*}\frac{d^{2}}{d\rho^{2}}_{\rho=0}(\delta p_{z})&=&\frac{d}{d\rho}_{|\rho=0}(p_{z}(x,z(x)+\xi)\frac{d\delta}{d\rho})=p_{z}(x,z( x))\frac{d^{2}\delta}{d\rho^{2}}_{\rho=0}=\frac{p_{z}}{qz}.\end{eqnarray*}

Before to consider the terms in the right-hand side observe that from \eqref{1411} it comes 
$$\delta'(x)=\frac{\rho^{2}}{2}\bigg(-\frac{z'}{qz^{2}}-\frac{q_{x} +q_{z}z'}{q^{2}z}\bigg)+o(\rho^{2}).$$

Hence $$\frac{d^{2}}{d\rho^{2}}_{\rho=0}(q\delta')=-\frac{p}{qz^{2}}-\frac {q_{x}+q_{z}\frac{p}{q}}{qz},$$
where we have replaced $z'$ by $p/q$, and $$\frac{d^{2}}{d\rho^{2}}_{\rho=0}(\delta\delta'q_{z})=0.$$
Finally

$$\frac{d^{2}}{d\rho^{2}}_{\rho=0}(\frac{\delta pq_{z}}{q})=
\frac{d}{d\rho}_{|\rho=0}(\frac{p(x,z(x))q_{z}(x,z(x)+\eta)d\delta/d\rho}{q(x,z(x))})=
\frac{pq_{z} }{q^{2}z}.$$

In conclusion we have obtained 
$$
\frac{p_{z}}{qz}=-\frac{p}{qz^{2}}-\frac{q_{x}}{qz}-\frac{pq_{z}}{q^{ 2}z}+\frac{pq_{z}}{q^{2}z}=-\frac{p}{qz^{2}}-\frac{q_{x}}{qz},
$$ that is to say $$p_{z}=-\frac{p}{z}-q_{x},$$ which, with d'Alembert's notation, is
$$B'=-\frac{p}{z}-A,$$
as we wanted to see.

\medskip

We now prove that $A'=B$.

Since the channel $FNM'K$ (Figure 17) is in equilibrium, the pressure along $NM'$ plus the pressure  following $FN$ must be  equal to the pressure along $FK$ plus the pressure  following $KM'$.

The pressure along $NM'$, according to \eqref{1611},  is equal to minus the integral of the second component of the acceleration, that by \eqref{2605b} is equal to
$$-\int_{0}^{\delta(x)}a^{2}(pB'+A'q)(x,z(x)+s)ds.$$

The pressure along $FN$ is, by \eqref{2905}, equal to $$\frac{a^{2}-U^{2}} {2}$$
where $U$ is the modulus of the velocity in $N$, that is  $U^{2}=a^{2}(p^{2}+q^{2})$.

The pressure along $FK$ is zero because the velocity far from the obstacle  is constant and therefore the acceleration in $FK$ is zero. 

Finally, the pressure along $KM'$ is equal to $$\frac{a^{2}-U'^{2}}{2}$$ where $U'$ is the modulus of the velocity in $M'$, that is    $U'^{2}=a^{2}(p'^{2}+q'^{2})=a^{2}( p^{2}+q^{2})+a^{2}\delta(x)\frac{\partial (p^{2}+q^{2})}{\partial z}$,
with this partial derivative valued at an intermediate point between $N$ and $M'$.

Equivalently 
\begin{eqnarray}\label{1511}U'^{2}-U^{2}=a^{2}\delta(x)\frac{\partial (p^{2}+q^{2})}{\partial z}_ {| (x,z(x)+\eta)},\qquad 0<\eta<\delta(x).\end{eqnarray}

Adding up the pressure on the channels and  using the mean value theorem for integrals,  we have 
\begin{eqnarray}\label{2605c}
&&-\int_{0}^{\delta(x)}a^{2}(pB'+A'q)(x,z(x)+s)ds+\frac{a^{2}-U^{2}} {2}\nonumber\\&=&
-\delta(x)a^{2} (pB'+A'q)(x,z(x)+\xi)+\frac{a^{2}-U^{2}}{2}\nonumber\\&= &\frac{a^{2}-U'^{2}}{2},\end{eqnarray}
with $0<\xi<\delta(x)$.  From \eqref{1511}  and \eqref{2605c}
we obtain
$$(pB'+A'q)(x,z+\xi)=\frac{1}{2}\frac{\partial (p^{2}+q^{2})}{\partial z}_ {| (x,z(x)+\eta)}.$$
Passing to the limit when $\delta(x)$ tends to zero, gives,  at the point $(x,z(x))$ (recall $B=q_{z}$, $B'=p_{z}$),
$$pB'+A'q=\frac{1}{2}\frac{\partial (p^{2}+q^{2})}{\partial z}=pB'+qB$$
and hence $A'=B$ as it should be proved.

\medskip

\subsection*{Exact differentials and Cauchy-Riemann equations}
As we said on page \pageref{2010}, d'Alembert intends 
to find functions $M(x,z), N(x,z)$ such that the differentials $M\,dx+N\,dz$ and $N\,dx-M\,dz$
are exact. The exactness of these differentials  means, in today's language, that the complex function $M-iN$ satisfies the Cauchy-Riemann equations
$$\frac{\partial M}{\partial x}=- \frac{\partial N}{\partial y},\qquad \frac{\partial M}{\partial y}=\frac{\partial N }{\partial x}$$ and it  is therefore holomorphic. Consequently $M-iN$ is a function only of $z$ and $M+iN$ is a function of $\bar{z}$. It follows that $M$ and $N$ are the sum and the difference of a function of $z$ and a function of $\bar{z}$, respectively. To reach this conclusion, d'Alembert's argument, \label{0507}which Euler later quotes in the article \cite{Euler3}, is as follows:

He notes that the forms $(M+iN)(dx-idy)$ and $(M-iN)(dx+idy)$ are exact, in particular $(M+iN)d\bar{z}=\frac{ \partial f}{\partial z}dz+\frac{\partial f}{\partial \bar{z}}d\bar{z}$
for a certain function  $f$, which satisfies therefore $\partial f/\partial z=0$ and hence it is a function  of $\bar{z}$ only. Thus  $M+iN$ is a function only of $\bar{z}$, say $\varphi(\bar{z})$. Similarly $M-iN$ is a function only of $z$, say  $\psi(z)$. Consequently \begin{eqnarray*}
M&=&\frac{1}{2}(\varphi(\bar{z})+\psi(z))\\
N&=&\frac{1}{2i}(\varphi(\bar{z})-\psi(z)).
\end{eqnarray*}

Reciprocally, it is easy to see that given two arbitrary functions $\varphi(\bar{z})$ and $\psi(z)$,
taking $M=\varphi(\bar{z})+\psi(z)$ and $N=-i (\varphi(\bar{z})-\psi(z))$, then the forms $Mdx +Ndy, Ndx-Mdy$ are exact.

%
%

\subsection{Euler's work on fluid theory}

Euler devotes three papers (\cite{Euler1}, \cite{Euler2} and \cite{Euler3}) to the theory of fluids in a spirit similar to that of d'Alembert. In the second of these articles, the famous Euler equations that regulate the movement of a fluid appear. For the particular case of homogeneous and stationary fluids, the first equation says that the form $udx+vdy+wdz$ is exact, where $u,v,w$ are the velocity components of the fluid and the second equation states that the divergence of the velocity field is zero.

Euler's argument to prove the second law is based on the principle of conservation of
mass. In a more current language it would essentially be the following\footnote{Nowadays it is usually proved using the divergence theorem that was stablished after Euler.}:

Consider an incompressible fluid of constant density one, occupying a region $C$ of the space. Denote by $X(x,y,z,t)$ the field of velocities
and let $(u,v,w)$ be its components.

To apply the law of conservation of mass that says us  that the rate of change of volume is zero,  we will calculate explicitly the volume of the region $C_{t}$ that the fluid fills after  a time $t$.

Denote by $\phi(x,y,z,t)$ the orbit of the point $(x,y,z)$, which satisfies
$$\frac{d\phi(x,y,z,t)}{dt}_{|t=0}=X(x,y,z,0),\quad \phi(x,y,z ,0)=(x,y,z).$$
Then the boundary of the region $C_{t}$ can be parameterized by $\phi(x,y,z,t)$ where $(x,y,z)$ is a point on the  boundary of $C$.

The linear part of the Taylor expansion of $\phi(x,y,z,t)$ at $t=0$ is
$$\phi(x,y,z,t)=(x,y,z)+t X(x,y,z,0).$$
One  can avoid the terms of higher order because in the subsequent calculations we will differentiate  with respect to $t$ at $t=0$.

To calculate the volume between the regions $C_{t}$ and $C$ one need to integrate the Jacobian of $\phi(x,y,z,t)$ with respect to $x,y,z$ which is

$$J_{\phi}=\left|\begin{array}{ccc}1+tu_{x} & tv_{x} & tw_{x} \\tu_{y} & 1+tv_{y} & tw_ {y} \\tu_{z} & tv_{z} & 1+tw_{z}\end{array}\right|=(1+tu_{x})(1+tv_{y})(1+tw_ {z})+O(t^{2})$$
and so
$$\frac{dJ_{\phi}}{dt}_{|t=0}=u_{x}+v_{y}+w_{z}=\operatorname{div}X(x,y,z, 0).$$

Denoting by $V_{t}$ the volume of $C_{t}\setminus C$ one has 
$$V_{t}=\int J_{\phi}(x,y,z,t)dx\,dy\,dz$$
and therefore
$$\frac{dV_{t}}{dt}_{|t=0}=\int \frac{dJ_{\phi}(x,y,z,t)}{dt}_{|t=0 }dx\,dy\,dz=\int \operatorname{div}X(x,y,z,0).$$
Since this derivative is the rate of change of the  volume of the fluid  at the boundary of $C$, it must be zero, and since this occurs  for any region $C$, the divergence of the velocity field must also be zero, that is, $u_{x}+v_{y}+w_{z}=0$, which is Euler's second equation.

\bigskip
In sections $LXIX$ and $LXX$ of \cite{Euler3}, Euler applies his equations to the case of a stationary fluid with constant density moving in planes perpendicular to a fixed axis. Choosing an orthogonal system of coordinates $(x,y,z)$ with $z$ in the direction of the fixed axis, the third component of the velocity will be zero.  In this situation the first Euler's equation reduces to the fact that
  the form $u\,dx+v\,dy$, where $u,v$ are the other two components of the velocity, must be exact, and the second one says that  $\frac{\partial u}{\partial x}+\frac{\partial v}{\partial y}=0$. This last equation is satisfied, by the equality of the cross derivatives, if the form $u\,dy-v\,dx$ is exact. The fact that the two forms $u\,dx+v\,dy$ and $u\,dy-v\,dx$ are exact, which is the same as saying that the Cauchy-Riemann equations are satisfied for the function $u-iv$, allows Euler to see, following d'Alembert's ``fort ingénieuse'' method (which we have reproduced on page \pageref{0507}), that $u,v$ must be the real part and the imaginary part of a function of $x+iy$,  as corresponds to the fact that the function $u-iv$ is holomorphic.

\medskip
A few years later, in a paper from 1777, \cite{Euler}, published in 1797, Euler starts checking   the validity of the Cauchy-Riemann equations for the complex derivative of a function. We outline his argument below.

Suppose that the function $f(z)$ has a primitive $F(z)$
such that $dF(z)/dz=f(z)$. If one  writes $f(z)=f(x+iy)=M(x,y)+iN(x,y)$ and
$F(z)=F(x+iy)=P(x,y)+iQ(x,y)$, one has 
\begin{eqnarray*}
dF&=&dP+idQ=f(z)dz=(M+iN)(dx+idy)\\&=&(Mdx-Ndy)+i(Ndx+Mdy),
\end{eqnarray*}
and therefore,
\begin{eqnarray*}
dP&=&Mdx-Ndy\\
dQ&=&Ndx+Mdy.
\end{eqnarray*}
By Clairaut's theorem of cross derivatives
\footnote{Euler does not explicitly cite Clairaut, but refers to a general property of functions that define an integrable form.}, it turns out $$\frac{\partial M}{\partial y}=- \frac{\partial N}{\partial x}, \qquad \frac{\partial N}{\partial y}= \frac{\partial M}{\partial x}.$$

Some years later, Cauchy in \cite{1814} also checks the Cauchy-Riemann equations for the derivative of a complex function,
as we discussed on page \pageref{0910}.

\subsection*{Relation to d'Alembert's equations}\label{1207}
We now apply the second Euler  equation in the case where a fluid meets an obstacle formed by a body of revolution, whose axis is the direction of the fluid, a situation considered by d'Alembert in \cite{dalembert}, section \S 44.

Suppose that the boundary of the revolution object is obtained by rotating around the $x$ axis the graph of a function in the $x,z$ plane, $(x,z(x))$\footnote{We keep d'Alembert's notation (see Figure 17).}. Let's also assume that the velocity components are $(q(x,z), p(x,z))$. After a rotation of  angle  $\theta$ the velocity at a point $(x,z,y)$ has components
\begin{eqnarray*}
q(x,z,y)&=&q(x,z)\\
p(x,z,y)&=&p(x,z)\cos \theta\\
r(x,z,y)&=&p(x,z)\sin\theta.
\end{eqnarray*}
  Therefore, the velocity as a function of the three variables $(x,z,y)$
is given by

\begin{eqnarray*}
q(x,z,y)&=&q(x,z)\\
p(x,z,y)&=&p(x,z)\frac{z}{\sqrt{y^{2}+z^{2}}}\\
r(x,z,y)&=&p(x,z)\frac{y}{\sqrt{y^{2}+z^{2}}}.
\end{eqnarray*}

So
\begin{eqnarray*}
\frac{\partial q}{\partial x}(x,z,y)&=&\frac{\partial q}{\partial x}(x,z)\\
\frac{\partial p}{\partial z}(x,z,y)&=&\frac{\partial p}{\partial z}(x,z)\cdot \frac{z}{\sqrt{ y^{2}+z^{2}}}+p(x,z)\cdot\frac{y^{2}}{(y^{2}+z^{2})^{3/2 }}\\
\frac{\partial r}{\partial y}(x,z,y)&=&p(x,z)\cdot \frac{z^{2}}{(y^{2}+z^{2) })^{3/2}}.
\end{eqnarray*}
Imposing now zero divergence and simplifying the notation one gets 
$$q_{x}+p_{z}\cdot \frac{z}{\sqrt{y^{2}+z^{2}}}+\frac{p}{\sqrt{y^{2} +z^{2}}}=0$$
and this equation restricted to the plane $y=0$
is
$$q_{x}=-p_{z}-\frac{p}{z},$$
which is the equation obtained by d'Alembert in \cite{dalembert}, and that we have  considered in Section \ref{1206}, equations \eqref{1206b}.

\medskip
When a fluid moving in planes perpendicular to a fixed axis, say the $y$ axis, 
encounters 
a cylinder whit 
axis in this direction, a situation considered by d'Alembert in section \S 73 of \cite{dalembert},
the velocity field can be written as $v=(q,0,p)$ with respect to an orthogonal coordinate
system $(x,y,z)$. Then Euler's equation
of null divergence gives $$q_{x}+p_{z}=0$$
which is equation \eqref{1206b} without the term $-p/z$. \label{3110}As we have commented in page \pageref{1207b} this result was announced without proof by d'Alembert.

 \section{Cauchy's theorem}\label{0410}
After this  review of the history of the Cauchy-Riemann equations we return to the fundamental problem about  the independence with respect to  the path of the integral of a complex function.
  
  As we explained at the end of Section \ref{0410b}, Cauchy deals with this topic in the article \cite{1825b} from the year 1825. We will try to understand how Cauchy reasoned in section \S 2 of this work. Unlike when he studied the case of a rectangle in \cite{1814} he will now not use the Cauchy-Riemann equations.

In fact, Cauchy's goal in \cite{1825b} is basically to extend the concept of integral between real limits to integrals between imaginary limits, that is to give a meaning to the expression
\begin{eqnarray}\label{1703}\int_{x_{0}+y_{0}i}^{X+Yi}f(z)dz.\end{eqnarray}
For this he considers the limits of the sums
$$S=\sum_{k}f(z_{k})(z_{k+1}-z_{k})$$
with $z_{k}=x_{k}+y_{k}i$, where the $x_{k}$ are a partition of $[x_{0},X]$ and the $y_{k}$ a partition of $[y_{0}, Y]$, when the partitions become more and  more fine. Cauchy remarks that there is no need for a single limit of these sums to exist.

One way to determine the previous partitions is through a continuous curve
$$x=\varphi(t), \,y=\chi(t),\quad t_{0}\leq t\leq T$$
with $\varphi(t)$ and $\chi(t)$ increasing functions, so that
\begin{eqnarray*}
\varphi(t_{0})&=&x_{0},\; \chi(t_{0})=y_{0}\\
\varphi(T)&=&X, \; \chi(T)=Y
\end{eqnarray*}
and taking then $x_{k}=\varphi(t_{k})$, $y_{k}=\chi(t_{k})$ where $t_{k}$ is a partition of $[t_{0} , T]$,

Cauchy, using an approximate version (as he says, {\em à trés-peu prés}) of the mean value theorem, that he states in the form  \begin{eqnarray*}
\varphi(t_{k+1})-\varphi(t_{k})&=&\varphi'(t_{k})(t_{k+1}-t_{k})\\
\chi(t_{k+1})-\chi(t_{k})&=&\chi'(t_{k})(t_{k+1}-t_{k}),
\end{eqnarray*}
writes the sum $S$ as
$$S=\sum_{k}f(\varphi(t_{k})+\chi(t_{k})i)(\varphi'(t_{k})+\chi'(t_{k}) i)(t_{k+1}-t_{k})$$
which tends, when refining the partition,  to the integral with real limits
\begin{eqnarray}\label{1703b}I=\int_{t_{0}}^{T}f(\gamma(t))\gamma'(t)dt, \quad \mbox{where } \gamma= \varphi + \chi i.\end{eqnarray}

In order that \eqref{1703b} provides a good definition of the integral \eqref{1703} it should be seen that it does not depend on the chosen curve  $\gamma$. It is now when Cauchy enunciates his famous theorem saying: if $f(x+yi)$ is
continuous\footnote{Cauchy says continuous but think of differentiable functions with a continuous derivative.} in the rectangle $[x_{0},X]\times [y_{0},Y]$ then the integral \eqref {1703b} is independent of the nature of the functions  $\varphi$ and $\chi$.

To prove this statement he will see that small variations of the functions $\varphi,\chi$ do not change  the value of the integral. Specifically, he considers two functions $u(t),v(t)$ vanishing  at $t_{0}$ and $T$ and  the curve $$\varphi(t)+\epsilon u( t)+(\chi(t)+\epsilon v(t))i=\gamma(t)+\epsilon W(t),$$
with $W(t)= u(t)+v(t)i$. Then he studies the increment $I_{\epsilon}-I$, where
$$I_{\epsilon}=\int_{t_{0}}^{T}f(\gamma(t)+\epsilon W(t))(\gamma'(t)+\epsilon W'(t) )\,dt,$$
   for sufficiently small values of $\epsilon$ so that the new curve does not leave the rectangle
   $[x_{0},X]\times [y_{0},Y]$, and he aims to see that the increment $I_{\epsilon}-I$ is zero.
  
  Cauchy affirms  that this increment  can be developed in series following the increasing powers of $\epsilon$ and says that the coefficient of $\epsilon$ is

  $$\int_{t_{0}}^{T}\bigg(f'(\gamma(t))\gamma'(t)W(t)dt+f(\gamma(t))W'(t )\bigg)\,dt$$
  which, integrating by parts, is seen to be  zero.

One way to find this coefficient of $\epsilon$ could be by calculating
$$\frac{dI_{\epsilon}}{d\epsilon}_{|\epsilon=0}.$$
Since Cauchy does not make it explicit, he could also be making a similar argument to the one he has used before, which consists of using the approximate equality
$$f(\gamma+\epsilon W)-f(\gamma)=\epsilon W f'(\gamma).$$
Then it would follow
\begin{eqnarray*}I_{\epsilon}-I&=&\int_{t_{0}}^{T}\bigg(f(\gamma+\epsilon W)\gamma'+ f(\gamma+\epsilon W) \epsilon W'\bigg)\,dt\\&=&\epsilon \int_{t_{0}}^{T}[W f'(\gamma)\gamma'+(f(\gamma)+\epsilon W f'(\gamma))W'] dt
\end{eqnarray*}
   which gives rise to the coefficient of $\epsilon$ considered before.

As a consequence of the vanishing of the coefficient of $\epsilon$,  Cauchy says that the increase of the integral $I_{\epsilon}-I$
    will reduce to an infinitely small of the second or higher order. Next he says: {\em Il est aisé d'en conclure que, si chacune des fonctions $x,y$, re\c coit successivement des accroissements infiniment petits du premier ordre dont la somme presente un accroissement fini, l'accroissement correspondant de $A+B\sqrt{-1}$ sera infiniment petit du premier ordre, c'est-à-dire nul}.
  
It is difficult to understand this argument of Cauchy. In fact, in none of the studies published on Cauchy's work  one can find a satisfactory explanation that clarifies it. A possible interpretation is that he is calculating the coefficients of $\epsilon^{2}, \epsilon^{3},$ etc. similarly to what he just did for the coefficient of $\epsilon$, either by directly calculating the successive derivatives of $I_{\epsilon}$ at $\epsilon=0$, or by considering the finite increments of successive order.
  
But another interpretation, perhaps more plausible, is that he is calculating the derivative of $I_{\epsilon}$ at successive points $\epsilon_{1}, \epsilon_{2},..$ through the increments
  
   $$\frac{I_{\epsilon_{i}+\epsilon}-I_{\epsilon_{i}}}\epsilon$$
   and in this way he would see that the derivative is zero at all these points and therefore that it is zero for all sufficiently small $\epsilon $. This would mean that the value of $I_{\epsilon}$ would not depend on $\epsilon$,  which is essentially what Cauchy wanted   to see.
  
   It is surprising that Cauchy does not realize, or at least does not comment, that his argument to see that $dI_{\epsilon}/d\epsilon$ is zero, at $\epsilon=0$, already gives that this derivative is zero for all sufficiently small $\epsilon $.
  
  Indeed, assuming like Cauchy that $f'$ is continuous we will have
  \begin{eqnarray*}
  &&\frac{d}{d\epsilon}\int_{t_{0}}^{T}f(\gamma+\epsilon W)(\gamma'+\epsilon W')\,dt=\\&&
  \int_{t_{0}}^{T}[f'(\gamma+\epsilon W)W(\gamma'+\epsilon W')+W'f(\gamma+\epsilon W)]\,dt=\\&&
  \int_{t_{0}}^{T}(f(\gamma+\epsilon W)W)'\,dt=0
  \end{eqnarray*}
  which directly proves that $I_{\epsilon}$ does not depend on $\epsilon$.

This reasoning, which also appears in Falk \cite{falk}, proves the following result: if $D$ is a convex domain of the plane and $f$ a function with continuous derivative in $D$ and $\gamma(t )$ and $\tilde{\gamma}(t)$ are two curves in $D$ that start and end at the same points, that is, they are defined in the interval $[t_{0},T]$ with $\gamma(t_{0})=\tilde{\gamma}(t_{0})$ and $\gamma(T)=\tilde{\gamma}(T)$
   then,
  $$\int_{t_{0}}^{T}f(\gamma(t))\gamma'(t)dt=\int_{t_{0}}^{T}f(\tilde{\gamma} (t)\tilde{\gamma}'(t)dt.$$

Indeed, it is only necessary to apply the previous argument with $W(t)=\tilde{\gamma}(t)-\gamma(t)$
and varying $\epsilon$ between $0$ and $1$. That this is really what Cauchy proves, namely the  path invariance of an integral for convex domains, is also discussed in \cite{bak}.
  
   \medskip
  On the other hand if we write $\gamma(t,\epsilon)=\gamma(t)+\epsilon(\tilde{\gamma}(t)-\gamma(t))$ 
  we have a homotopy between $\gamma$
and $\tilde{\gamma}$, and the integral along the curves $\gamma(\cdot,\epsilon)$ does not depend on $\epsilon$.

This suggests that the independence of an integral with respect to the path of integration will hold for any homotopy.
Specifically, it can be shown with arguments similar to the previous ones, that if we have two homotopic curves, with the same beginning points  and ending points, in a domain $D$ where $f$
has a continuous derivative, then the integral of $f$ along the two curves has the same value. This is done in \cite{ber} when the homotopy is of class ${\cal C}^{2}$ and in \cite{flato} when the homotopy is only continuous
   and without assuming that $f'$ is continuous.

  In fact it seems that Cauchy does not find his reasoning clear enough  since he immediately reformulates it in terms of the calculus of variations. The total variation of the integral $I$
   when $\gamma$ receives an increment $\delta \gamma=\epsilon W$
is
\begin{eqnarray*}\delta I&=&\int_{t_{0}}^{T}\delta(f(\gamma)\gamma')\,dt=\int_{t_{0}}^{T }(f'(\gamma) \delta \gamma\cdot \gamma' +f(\gamma)\delta\gamma')\,dt\\&=&\epsilon \int_{t_{0}}^{T } (f'(\gamma)\gamma' W+f(\gamma)W')dt\end{eqnarray*}
which is zero, as we have seen. Here it has been used that $\delta$ behaves like a derivative, that is  $(\delta \gamma)'=\delta\gamma'$ and $\delta(f(\gamma))=f' (\gamma)\delta\gamma$.

Later, Casorati in 1868 also justified Cauchy's theorem with this variational argument (\cite{casorati}, p.365).


\subsection*{From Cauchy to Goursat}
We note that Cauchy always assumed that the functions he considered had a continuous derivative. In fact he was not even talking about the derivative and was simply referring to complex functions. But to reach the conclusion of Cauchy's theorem it is not necessary to assume that the considered function has a continuous derivative, but only that it is holomorphic, as it was proved by Goursat in 1900.

In 1875, one year after the publication of Cauchy's article \cite{1825a}, Briot and Bouquet in \cite{BB}
published a proof of Cauchy's theorem, for holomorphic functions, without assuming that they have a continuous derivative. They first do this for a star shaped domain and then they go  to the general case by subdividing the initial region into smaller ones. However, their argument in the star shaped case is not entirely correct because they assumed that there was uniformity in the definition of the derivative (see \cite{bak}). 

In 1884 Goursat in \cite{goursat1} presents a rigorous demonstration of Cauchy's theorem that he considers simpler than the previous ones, but still assuming that the holomorphic function  has a continuous derivative. Finally in 1900, at  Osgood's request, Goursat himself in \cite{goursat2} complements the proof of 1884 so that it is valid without the condition of the continuity of the derivative. It can be regarded as the first correct proof of Cauchy's theorem assuming only that the function is holomorphic.
 
  As Goursat comments, this result gives rise to Cauchy's integral formula for a holomorphic function and from here the theory of analytic functions can be developed from Cauchy's point of view, in particular showing that a holomorphic function already has continuous derivative.
  
  \section{On Green's formula}\label{1110}
  As we said in the Introduction, Green's formula, which plays a fundamental role in the theory of functions of a complex variable, together with the Cauchy-Riemann equations is now the usual way to prove Cauchy's theorem. For this reason we will do a brief historical review of this formula.
 
  Green's formula states that
\begin{eqnarray}\label{0707b}\int_{\partial R} Pdx+Qdy=\int_{R} (\frac{\partial Q}{\partial dx}-\frac{\partial P}{dy })dx\,dy\end{eqnarray}
where $P(x,y), Q(x,y)$ are regular functions of $x,y$ and $R$ is a region of the plane.

The reasoning that allows  to obtain Cauchy's theorem
from this equality and the Cauchy-Riemann equations is the following:

If $f(z)=u(z)+iv(z)$ is a holomorphic function in the region $R$ then $u_{x}=v_{y}$ and $u_{y}=-v_{ x}$ and therefore integrating $f(z)$ along the closed path $\partial R$, and applying Green's formula it  results
\begin{eqnarray*}
\int_{\partial R} f(z)dz&=&\int_{\partial R} (u+iv)(dx+idy)=\int_{\partial R} (u \,dx-v\,dy) + i \int_{\partial R} (v\,dx+u\,dy) \\
&=&- \int_{R} (v_{x}+u_{y})\,dx\,dy+i \int_{R} (u_{x}-v_{y})\,dx\,dy =0.
\end{eqnarray*}


  Let's see below how Green's formula has appeared in the literature. We have already commented in the Introduction that the first time that it was stated in the form \eqref{0707b} is in Cauchy's article from the year 1846, \cite{cauchy1846}.

   The first published proof of Green's formula is due to Riemann in his 1851 dissertation, \cite{riemann}, which we reproduce in Section \ref{0407}. Later on this formula appears in Casorati \cite{casorati}, p.381.

\medskip

On the other hand, Green's formula is the version in the plane  of the divergence theorem in the space, that states 
$$ \int_{\partial S}\langle X,N\rangle dS=\int_{S}(\frac{\partial P}{\partial x}+\frac{\partial Q}{\partial y}+ \frac{\partial R}{\partial z})\,dx\,dy\,dz,$$
where $X=(P,Q,R)$ is a regular field in $\R^{3}$, $S$ a region of the space and $N$ the unit normal exterior vector
  to the boundary of $S$. Indeed, applying this formula to the field $(Q,-P,0)$ and to a region formed by the points $(x,y,t)$ with $(x,y)\in S\subset\R^{2 }$ and $0\leq t\leq \epsilon$ we get Green's formula.

  The divergence theorem was proved in a particular case by Gauss in 1813, \cite{gauss}. After that, in 1828, Green
  in  \cite{green} proved it for a special type of vector fields. In its general version it was established by Ostrogadsky in 1831.

  \medskip
Likewise, Green's formula is also the version for the plane of the rotational theorem,
that stays 
$$\int_{S} \langle (\frac{\partial R}{\partial y}- \frac{\partial Q}{\partial z},\frac{\partial P}{\partial z}- \frac{\partial R}{\partial x},\frac{\partial Q}{\partial x}- \frac{\partial P}{\partial y} ), N\rangle dS=\int_{\partial S }\langle (P,Q,R), T\rangle ds$$
where $P,Q,R$ are regular functions, $S$ a surface of the space with unit normal vector $N$ and $T$ the unit tangent vector to the boundary of $S$.

Historically this theorem, also known as Stokes' theorem because Stokes  himself proposed it in an exam (see \cite{katz}), was not proved until 1861 by Hankel (\cite{hankel}, p.34 ) who based his proof precisely on Green's formula, for which he quotes Riemann \cite{riemann2}.

 \subsection{Gauss's formula for the volume of a body}
Gauss in the 1813 paper \cite{gauss1} proves a formula for the volume of a body that can be considered as a particular case of the divergence theorem.

Essentially, what Gauss does can be explained as follows. When a surface $S$ is parameterized as the graph of a function of the form $\varphi(x,y)=(x,y,z(x,y))$
then the ratio between the area element $dS$ of the surface and the area element $dx\,dy$ in the plane  is $$\pm \cos\theta\, dS=dx \,dy$$ where $\theta$ is the angle between the normal vector to the surface, $(-z_{x},-z_{y},1)$ or $(z_{x},z_{y},-1 )$, and the positive part of the $z$ axis. The sign $+$ corresponds to the points where $\theta$ is acute and the sign $-$ when $\theta$ is obtuse, depending on the normal vector we have chosen.

Gauss gives the above relation between $dS$ and $dx\,dy$ without further comment, but the justification is clear if we think that a small rectangle of sides $\Delta x, \Delta y$ centered at a point $P_{0}$ in the plane $z=0$, is the projection of a rectangle in the tangent plane to the surface at the point $P$ that is projected onto $P_{0}$. The ratio between the area $\Delta x\cdot\Delta y$ of the rectangle in the plane $z=0$ and the area $\Delta S$ of the rectangle in the tangent plane is $\pm \cos\theta \Delta S=\Delta x\cdot \Delta y$, where $\theta$ is the angle formed by these two planes.
  
Gauss considers a compact surface that encloses a region $K$ of the space. This allows him to choose a normal exterior vector  to the surface at each point. He observes that each point $P_{0}$ of the plane $z=0$ is the projection of an even number of points in the surface, and that the sign of the cosine of the angle that forms the normal exterior vector with the axis $z$ alternates from one point to another.

\begin{figure}[h]
\begin{center}
\includegraphics[width=.5\textwidth]{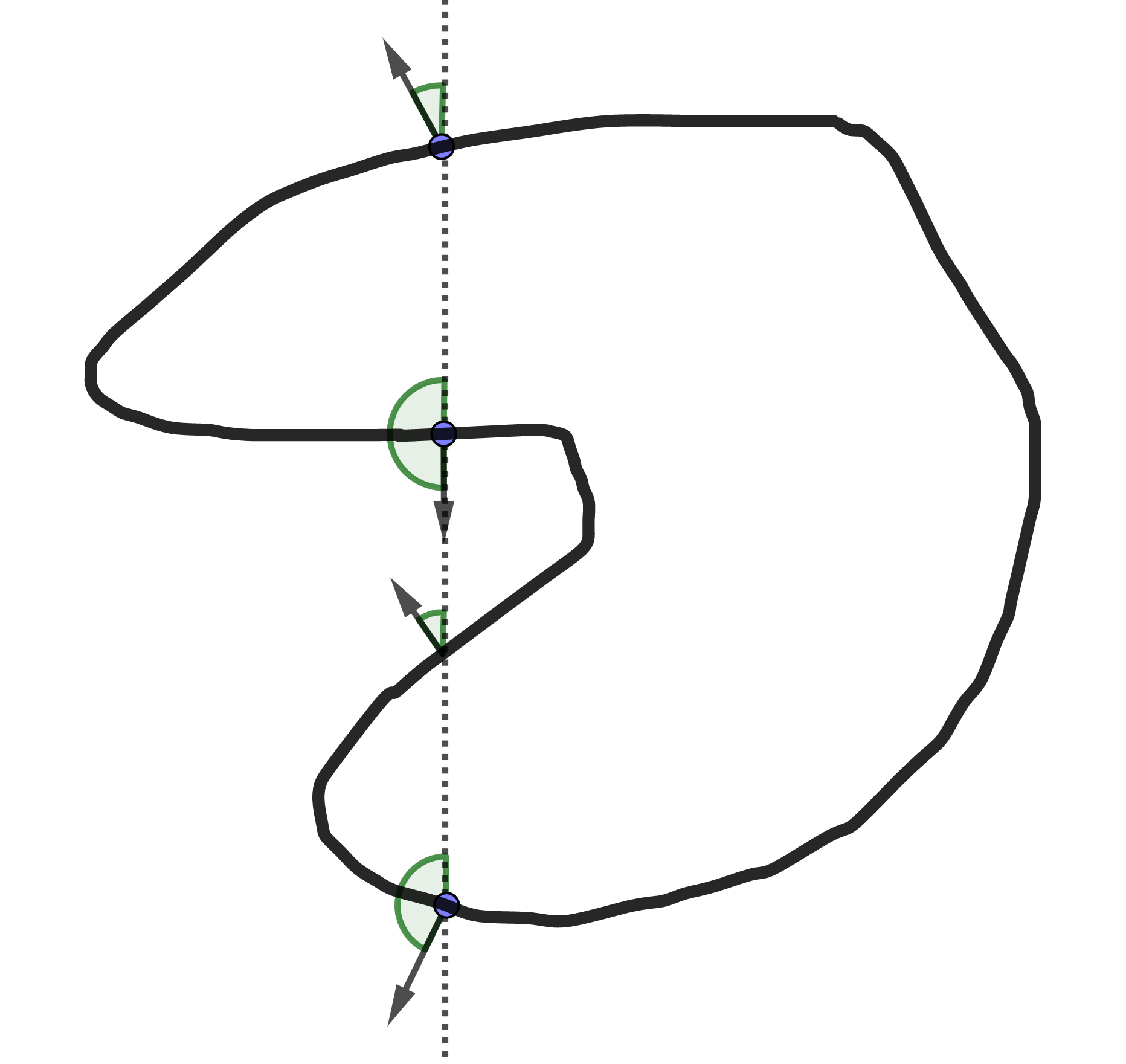}
\end{center}
\end{figure}

For calculating the volume $V$ of the body $K$ we can assume that it is  decomposed into areas limited by graphs of functions.
In fact, the problem reduces to calculate the volume of a body bounded by the graph of two functions $z_{1}(x,y), z_{2}(x,y)$, defined on the same region $D $ of the plane $z=0$, with $z_{1}(x,y)\geq z_{2}(x,y)$. Let $S_{1}$ be the surface given by the graph of $z_{1}$ and $S_{2}$ the one given by the graph of $z_{2}$, so that the total surface of the body is $S=S_{1}\cup S_{2}$. The normal exterior vector makes an acute angle with the positive part of the $z$ axis at the points of $S_{1}$ and an obtuse one at the points of $S_{2}$.

So \begin{eqnarray*}V&=&\int_{D}z_{1}(x,y)dx\,dy-\int_{D}z_{2}(x,y)dx\,dy\\& =&\int_{S_{1}}z_{1}\cos\theta\,dS_{1}+\int_{S_{2}}z_{2}\cos\theta\,dS_{2}=\int_ {S}z\cos\theta dS.\end{eqnarray*}

Analogously, changing  the plane $z=0$ by the planes $x=0$ and $y=0$ 
one would  obtain, respectively,  $$V=\int_{S}x \cos\theta dS=\int_{S}y \cos\theta dS, $$ where $\theta$
is the angle of the normal with the corresponding axes.
Gauss gives these three equalities on page 6 of \cite{gauss1}. 

Note that they  are a particular case of the divergence theorem applied to the fields $X=(x,0,0)$, $Y=(0,y,0)$ and $Z=(0,0,z )$, respectively. Indeed, the above equalities can be written as
\begin{eqnarray*}
V&=&\int_{K}\operatorname{div}X dx\,dy\,dz=\int_{S}\langle N,X\rangle dS\\
V&=&\int_{K}\operatorname{div}Y dx\,dy\,dz=\int_{S}\langle N,Y\rangle dS\\
V&=&\int_{K}\operatorname{div}Z dx\,dy\,dz=\int_{S}\langle N,Z\rangle dS
\end{eqnarray*}
where $N$ is the unit normal exterior vector to the surface.
\subsection{Ostrogradsky's divergence theorem}
With exactly the same kind of arguments as Gauss, Ostrogradsky in \cite{ostro} proves the general version of the divergence theorem. It is enough  to see that for a function $R=R(x,y,z)$ it holds
$$\int_{K}\frac{\partial R}{\partial z}dx\,dy\,dz=\int_{S}R\cos\theta dS,$$
where $\cos\theta$ is the third component of $N$,  the unit normal exterior  vector to $S$. This equality is deduced, as we discussed above, by reducing the problem to the case where $K$ is bounded by the graphs
of two functions $z_{1}(x,y), z_{2}(x,y)$. Then,

\begin{eqnarray*}\int_{K}\frac{\partial R}{\partial z}dx\,dy\,dz&=&\int_{D}(R(x,y,z_{1})- R(x,y,z_{2}))\,dx\,dy\\&=&\int_{S_{1}}R(x,y,z_{1})\cos\theta\,dS_{ 1}+\int_{S_{2}}R(x,y,z_{2})\cos\theta\,dS_{2}\\&=&\int_{S}R\cos\theta dS.\end{eqnarray*}

In conclusion, given a field $X=(P,Q,R)$ one gets
$$\int_{K}(\frac{\partial P}{\partial x}+\frac{\partial Q}{\partial y}+\frac{\partial R}{\partial z})\,dx \,dy\,dz= \int_{S}\langle X,N\rangle dS,$$
which is the statement of the divergence theorem.

This paper of Ostrogradsky that we are discussing was published in 1831, but read in 1828. Meanwhile, in 1829, Poisson, who had been Ostrogradsky's supervisor in Paris, published essentially the same proof in \cite{poisson}. To have more historical details about
the divergence theorem see \cite{katz}.

\subsection{The study of the electric field according to Green}
Green's formula \eqref{0707b} does not appear in Green's paper \cite{green}. In fact, in this article the author proves the formula
\begin{eqnarray}\label{2903}
\int dx\,dy\,dz\,U \delta V+\int d\sigma U(\frac{dV}{dw})=\int dx\,dy\,dz\,V \delta U+\int d \sigma V(\frac{dU}{dw})
\end{eqnarray}
where $U,V$ are differentiable functions of $x,y,z$ defined on a region of the space, where the  triple integrals are extended, $d\sigma$ is the area element of the boundary of this region, $dw$ is the length element in the direction of the interior normal, and $\delta$ represents the Laplacian. The above equality gives  a relationship between the density of an electric fluid on the surface of a body and the potentials inside and outside this surface.

We note that \eqref{2903} is a consequence of divergence theorem
\begin{eqnarray}\label{1204}\int \langle X, N \rangle\, d\sigma=\int \operatorname{div}X \,dx\,dy\,dz\end{eqnarray}
where $X$ is a field in the considered region and $N$ is the outer normal  vector to the boundary of this region. Indeed, since $$\frac{dU}{dw}=-\langle \nabla U, N \rangle,\quad \frac{dV}{dw}=-\langle \nabla V, N \rangle, $$
taking $X=U \cdot \nabla V$ and denoting now the Laplacian by $\Delta$,  one gets  $$\operatorname{div}X=\langle \nabla U, \nabla V\rangle+ U\cdot \Delta V.$$
It  turns out that
$$\int U\frac{dV}{dw}=-\int \langle U\cdot \nabla V,N\rangle\,d\sigma=-\int \langle \nabla U,\nabla V\rangle dx \,dy\,dz-\int U\cdot \Delta V\, dx\,dy\,dz$$
and, similarly,
$$\int V\frac{dU}{dw}=-\int \langle V\cdot \nabla U,N\rangle\,d\sigma=-\int \langle \nabla U,\nabla V\rangle dx \,dy\,dz-\int V\cdot \Delta U\, dx\,dy\,dz$$ and these two equalities give \eqref{2903}.

\medskip
For fields $X$ that  can be written as $X=U\cdot \nabla V-V\cdot \nabla U$, for certain functions $U,V$, the above calculation shows that \eqref{1204} is a consequence of \eqref{2903}.

 \section{Riemann's contributions}\label{0407}
   A few years later of the publication of the works of Cauchy and of Gauss and Ostrogradsky, that we have discussed, Riemann in his famous thesis on the foundations of the general theory of functions of a complex variable, presented in 1851, \cite{riemann}, rediscovers the Cauchy-Riemann equations and gives
   the first direct proof of Green's formula. Let us follow his arguments.

\subsubsection*{Cauchy-Riemann Equations}
Riemann defines a holomorphic function\footnote{Riemann does not use this word.} $w=w(z)$
when  the differential quotient $dw/dz$ is independent of the value of the differential $dz$.

This definition can be interpreted as follows: the image by $w$ of a curve $z(t)$ in the complex plane, with $z(0)=a$, verifies

\begin{eqnarray}\label{1406}\frac{dw}{dz}=\frac{dw/dt}{dz/dt}=\frac{(w\circ z)'(0)}{z'( 0)}\end{eqnarray} and this quotient does not depend on the curve $z(t)$ passing through $a$. Therefore, this ratio of tangent vectors considered  as complex numbers
 can be  called  $w'(a)$.

It is easily seen that this Riemann's definition coincides with the usual one, that is,

$$w'(a)=\lim_{z\to a}\frac{w(z)-w(a)}{z-a}.$$

From the relation \eqref{1406} it follows that the function $w(z)$ is conformal.
Indeed, if $\gamma(t), \sigma(t)$ are two curves with $\gamma(0)=\sigma(0)$
then
$$\frac{\gamma'(0)}{\sigma'(0)}=\frac{(w\circ \gamma)'(0)}{(w\circ \sigma)'(0)}$$
equality which means that the function $w$ preserves the angle between the tangents to the curves $\gamma$ and $\sigma$ at a common point.

Riemann then proves that the function $w$
satisfies the Cauchy-Riemann equations.
He deduces it as follows: putting $z=x+iy$ and $w=u+iv$ the differential quotient $dw/dz$ is written as
\begin{eqnarray*}
\frac{du+i\,dv}{dx+i \,dy}&=&\frac{(\frac{\partial u}{\partial x}+i\,\frac{\partial v}{\partial x})dx+(\frac{\partial v}{\partial y}-i\,\frac{\partial u}{\partial y})i\,dy}{dx+i\,dy}.
\end{eqnarray*}

Since this quotient does not depend on the value $dx+i\,dy$, considered as the tangent vector to a curve, and taking the curves $x=constant$ and $y=constant$,
it turns out that
$$\frac{\partial u}{\partial x}+i\,\frac{\partial v}{\partial x}=\frac{\partial v}{\partial y}-i\,\frac{ \partial u}{\partial y}$$
or, equivalently, 
$$\frac{\partial u}{\partial x}=\frac{\partial v}{\partial y}\qquad\mbox{and}\qquad \frac{\partial v}{\partial x}=- \frac{\partial u}{\partial y}.$$

\subsubsection*{Green's formula}
Regarding Green's formula, in section 7 of his thesis, Riemann proves that
\begin{eqnarray}\label{1606}\int \bigg(\frac{\partial X}{\partial x}+\frac{\partial Y}{\partial y}\bigg)dT=-\int (X \cos\xi+Y\cos\eta)ds\end{eqnarray}
where $X,Y$ are 
continuous functions\footnote{Riemann talks about continuous functions although it is clear that he needs more assumptions.} of $x,y$ defined on a plane  region $T$, and writes the unit interior normal vector at a point of the boundary of $T$ as $(\cos\xi, \cos\eta)$. It is understood that the integral on the left-hand side is extended to $T$, with $dT=dx dy$ being the area element, and the integral on the right-hand side is extended to the boundary of $T$ with length element $ ds$.

We note that the above equality,  which is the divergence theorem in the plane, gives rise to Green's formula applied to the field $(Y, -X)$.

The proof of \eqref{1606} reduces to show the  equality

\begin{eqnarray}\label{1906b}\int \frac{\partial X}{\partial x}dT=-\int X\cos\xi\, ds.\end{eqnarray}

For this purpose Riemann divides the region $T$ into trapezoidal regions obtained by cutting $T$ by lines parallel to the  $x$ axis, equally separated from each other, so that it is enough to prove \eqref{1906b} for everyone of these trapezoids.

Assume that $\Delta y$ is small enough so that the two boundaries 
  of these trapezoids can be parameterized
in the form $(x(y),y)$. We fix one of these trapezoids and parameterize the left border by $(x_{1}(y),y)$ and the right one by $(x_{2}(y),y)$.

The orientation of the boundary of $T$, that must be travelled leaving $T$ on the left side, implies that the interior normal vector $N$ is given by

\begin{eqnarray}\label{1906c}N=(1,-\frac{dx_{1}}{dy})=\|N\| (\cos\xi, \cos\eta)\end{eqnarray}
on the left border of the trapezoid, and by
\begin{eqnarray}\label{1906d}N=(-1,\frac{dx_{2}}{dy})=\|N\| (\cos\xi, \cos\eta)\end{eqnarray}
on the right border. In particular $\cos\xi$ is positive in the first case and negative in the second one.

If for a fixed trapezoid $y$ varies in a certain interval $[a,b]$ we take the partition determined for each $n$ by the points $y_{k}=a+k\frac{b-a}{n} $, $k=0,\dots, n$, so that $\Delta y=y_{k+1}-y_{k}$.

Then
\begin{eqnarray}\label{1906}
\int \frac{\partial X}{\partial x}dT&=&\int_{a}^{b}\int_{x_{1}(y)}^{x_{2}(y)}\frac{ \partial X}{\partial x}\,dx\,dy=
\int_{a}^{b}(X(x_{2}(y),y)-X(x_{1}(y),y))\,dy\nonumber\\
&=&\lim_{n\to \infty}\Delta y \sum_{k=0}^{n}(X(x_{2}(y_{k}),y_{k})-X(x_{ 1}(y_{k}),y_{k})).
\end{eqnarray}

On the other hand, in order to calculate $\int X\cos \xi\,ds$ on the boundary  of the trapezoid it is only necessary to consider the non-horizontal part of this boundary  since on the horizontal parts it is $\cos\xi=0$ . The partition we have considered given by the values $y_{k}$ determines a partition of the left part of the boundary  given by the points $(x_{1}(y_{k}),y_{k})$. Denoting by $s_{k}$ the length over this boundary from the point $(x_{1}(a),a)$ to the point $(x_{1}(y_{k},y_{k}))$, we have that $\int X\cos \xi\,ds$ extended to the left border  is
$$\lim_{n\to\infty}\sum_{k=0}^{n}(s_{k+1}-s_{k})X(x_{1}(y_{k}),y_{ k})\cos \xi(x_{1}(y_{k}),y_{k}).$$

Now, by \eqref{1906c}, we have
\begin{eqnarray*}s_{k+1}-s_{k}&=&\int_{y_{k}}^{y_{k+1}}\sqrt{1+(\frac{dx}{dy })^{2}}dy=\Delta y \bigg(\sqrt{1+(\frac{dx}{dy})^{2}}\bigg)_{y=\tau_{k}} \\&=& \frac{\Delta y}{(\cos\xi)_{y=\tau_{k}}}, \qquad y_{k}\leq \tau_{k}\leq y_{k+1 }.\end{eqnarray*}

Therefore the integral we are considering, on the left border, equals to 
$$\int X\cos\xi \,ds=\lim_{n\to\infty}\Delta y\sum_{k=0}^{n}X(x_{1}(y_{k}),y_ {k}).$$

Similarly, on the right border, using now \eqref{1906d}, it results
$$\int X\cos\xi \,ds=- \lim_{n\to\infty}\Delta y\sum_{k=0}^{n}X(x_{2}(y_{k}), y_{k}).$$

These two equalities added themselves,  together with \eqref{1906} prove
\eqref{1906b}.

Analogously it would be shown that $\int\frac{\partial Y}{\partial y}dT=-\int Y \cos\eta ds$ and we would obtain \eqref{1606}.

\subsection*{Chronological order of results}

\begin{itemize}

\item[\bf 1740]
\begin{itemize}\item Clairaut: {\em Sur l'intégration ou la construction des équations différentielles du premier ordre}, \cite{Clairaut1740}.\end{itemize}
\item[\bf 1743]
\begin{itemize}\item Clairaut: {\em Théorie de la figure de la Terre, tirée des principes de l'Hydrostatique}, \cite{Clairaut1743}.\end{itemize}

\item[\bf 1752]
\begin{itemize}\item d'Alembert: {\em Essai d'une nouvelle théorie de la résistance des fluides},  \cite{dalembert}.\end{itemize}

\item[\bf 1757]
\begin{itemize}\item Euler: {\em Equilibre et mouvement  des fluides I, II, III,} \cite{Euler1}, \cite{Euler2}, \cite{Euler3}.\end{itemize}

\item[\bf 1768]
\begin{itemize}\item d'Alembert: {\em Sur l'équilibre des fluides}, \cite{dalembert1768}.\end{itemize}

\item[\bf 1797]
\begin{itemize}\item Euler: {\em Ulterior Disquisitio de Formulis Integralibys Imagirariis.} Written in 1777, \cite{Euler}.\end{itemize}

\item[\bf 1811]
\begin{itemize}\item Gauss's letter to  Bessel, \cite{gauss}.\end{itemize}

\item[\bf 1813]
\begin{itemize}\item Gauss: {\em  Theoria attractionis corporum Sphaeroidicorum Ellipticorum Homogeneorum, Methodo nova tractata}, \cite{gauss1}.\end{itemize}

\item[\bf 1814]
\begin{itemize}\item Cauchy: {\em Mémoire sur les integrales définies}. Presented at the  Académie; it was not published until 1827  in \cite{1814}.\end{itemize}

\item[\bf 1822]
\begin{itemize}\item Cauchy: {\em Mémoires sur les intégrales définies oú l'on fixe le nombre et la nature des constantes arbitraires et des fonctions arbitraires}, \cite{1822b}.\end{itemize}

\item[\bf 1823]
\begin{itemize}\item Cauchy: {\em Résumé des le\c cons sur le calcul infinitésimal}, \cite{resume}.\end{itemize}

\item[\bf 1824]
\begin{itemize}\item Cauchy: {\em Mémoire sur les integrales définies prises entre des limits imaginaires}. Presented at the  Académie; it was not published until  1825  (\cite{1825b}) and 1874 (\cite{1825a}).\end{itemize}

\item[\bf 1828]
\begin{itemize}\item Green: {\em An Essay on the Application of Mathematical Analysis to the theories of Electricity and Magnetisme}, \cite{green}.\end{itemize}

\item[\bf 1828]
\begin{itemize}\item Ostrogradsky: {\em Note sur la Théorie de la Chaleur}, \cite{ostro}.\end{itemize}

\item[\bf 1829]
\begin{itemize}\item Poisson: {\em Mémoire sur l'Équilibre et le Mouvement des Corps Élastiques}, \cite{poisson}.\end{itemize}

\item[\bf 1846]
\begin{itemize}\item Cauchy: {\em Sur les intégrales qui s'étendent à tous les points d'une courbe fermée},  \cite{cauchy1846}.\end{itemize}

\item[\bf 1851]
\begin{itemize}\item Riemann: {\em Grundlagen für eine allgemeine Theorie der
Functionen einer veränderlichen complexen
Grösse}, \cite{riemann}.  \end{itemize}

\item[\bf 1857]
\begin{itemize}\item Riemann: {\em Lehrsätze aus der analysis situs für die Theorie der Integrale von zweigliedrigen vollständigen Differentialien}, \cite{riemann2}.\end{itemize}

\item[\bf 1861]
\begin{itemize}\item Hankel: {\em Zur allgemeinen Theorie der Bewegung der Flüssigkeiten}, \cite{hankel}.\end{itemize}

\item[\bf 1868]
\begin{itemize}\item Casorati {\em Teorica delle funzioni di variabili complesse},  \cite{casorati}.  \end{itemize}

\item[\bf 1875]
\begin{itemize}\item Briot-Bouquet: {\em Théorie des fonctions elliptiques}, \cite{BB}. \end{itemize}

\item[\bf 1883]
\begin{itemize}\item Falk: {\em Extrait d'une lettre adressée à M. Hermite}, \cite{falk}. \end{itemize}

\item[\bf 1884]
\begin{itemize}\item Goursat: {\em D\'{e}monstration du th\'{e}or\'{e}me de {C}auchy},  \cite{goursat1}. \end{itemize}

\item[\bf 1900]
\begin{itemize}\item Goursat: {\em Sur la d\'{e}finition g\'{e}n\'{e}rale des fonctions analytiques, d'apr\`es
              {C}auchy},  \cite{goursat2}. \end{itemize}
\end{itemize}

 \bibliographystyle{amsplain}

\bibliography{BibCauchyLatin1}

{\sc Department of Mathematics, Universitat Autònoma de Barcelona, 08193, Bellaterra, Barcelona, Catalonia}

{\em E-mail address}: julia.cufi@uab.cat, agusti.reventos@uab.cat
\end{document}